\documentclass{ajour}

\usepackage{amssymb}

\begin{document}

%


\authorrunninghead{Petr Vojt\v echovsk\'y}
\titlerunninghead{Generators of Simple Moufang Loops}






\title{Generators of Nonassociative Simple Moufang Loops over Finite Prime
Fields}
\author{Petr Vojt\v echovsk\'y}
\affil{Department of Mathematics, Iowa State University, Ames, IA 50011,
U.S.A.} \email{petr@iastate.edu}

\abstract{ We present an elementary proof that the nonassociative simple
Moufang loops over finite prime fields are generated by three elements. In the
last section, we conclude that integral Cayley numbers of unit norm are
generated multiplicatively by three elements.}



\def\to{\longrightarrow}
\def\melm#1#2#3#4{\left(\begin{array}{cc}#1 & #2\\ #3 & #4\end{array}\right)}
\def\uelm#1{\melm{1}{#1}{0}{1}}
\def\lelm#1{\melm{1}{0}{#1}{1}}
\def\delm#1#2{\melm{#1}{0}{0}{#2}}
\def\aelm#1#2{\melm{0}{#1}{#2}{0}}
\def\det#1{\textrm{det}\, #1}
\def\char#1{\textrm{char}\, #1}
\def\image#1{\textrm{Im}\, #1}
\def\myspace{\;\;}

\begin{article}


\section{Simple Moufang Loops}

\noindent The first class of nonassociative simple Moufang loops was
discovered by L.~Paige in $1956$ \cite{Paige}, who investigated Zorn's and
Albert's construction of simple alternative rings. M.~Liebeck proved in $1987$
\cite{Liebeck} that there are no other finite nonassociative simple Moufang
loops. We can briefly describe the class as follows:

For every finite field $\mathbb F$, there is exactly one simple Moufang loop.
Recall Zorn's multiplication
\begin{displaymath}
    \melm{a}{\alpha}{\beta}{b}\melm{c}{\gamma}{\delta}{d}=
    \melm{ac+\alpha\cdot\delta}{a\gamma+\alpha d-\beta\times\delta}
    {\beta c+b\delta+\alpha\times\gamma}{\beta\cdot\gamma+bd},
\end{displaymath}
where $a$, $b$, $c$, $d\in \mathbb F$, $\alpha$, $\beta$, $\gamma$,
$\delta\in{\mathbb F}^3$, and where $\alpha\cdot\beta$ (resp.
$\alpha\times\beta$) denotes the dot product (resp. vector product) of $\alpha$
and $\beta$. Probably the easiest way to think of this multiplication is to
consider the usual matrix multiplication with the added antidiagonal matrix
\begin{displaymath}
    \melm{0}{-\beta\times\delta}{\alpha\times\gamma}{0}.
\end{displaymath}
We define the determinant of
\begin{displaymath}
    M=\melm{a}{\alpha}{\beta}{b}
\end{displaymath}
by $\det{M}=ab-\alpha\cdot\beta$. Then $\mathcal L=\{M;$ $\det{M}\ne 0\}$ turns
out to be a nonassociative Moufang loop, and so does $\mathcal M=\{M;$
$\det{M}=1\}$. One can show that $Z(\mathcal L)$, the center of $\mathcal L$,
consists of all elements
\begin{displaymath}
    \melm{a}{0}{0}{a},
\end{displaymath}
for $0\ne a\in \mathbb F$. Thus $Z(\mathcal M)$ is at most a two-element
subloop of $\mathcal M$. More precisely, $Z(\mathcal M)$ is trivial if and only
if $\char{\mathbb F}$, the characteristic of $\mathbb F$, equals $2$ (cf. Lemma
$3.2$ \cite{Paige}). Finally, $\mathcal M/Z(\mathcal M)$ was found to be simple
(and nonassociative) in \cite{Paige}.

Obviously, a finite simple Moufang loop is either associative (whence a finite
simple group), or nonassociative---an element of the class introduced above.
Each finite simple group is known to be generated by just two elements. (See
\cite{Wilson} for more details.) The proof of this fact depends heavily on the
classification of finite simple groups. As we have seen, nonassociative simple
Moufang loops admit a much simpler classification. They cannot be generated by
two elements only, since every Moufang loop is diassociative (cf.
\cite{Moufang}, \cite{Bruck}, or \cite{Pflugfelder}).

In this paper, we show that $\mathcal M/Z(\mathcal M)$ is generated by three
elements when $\mathbb F$ is a (finite) prime field. The approach taken here is
elementary. Possible generalizations of this result would probably require
more detailed methods than those used here, or a completely different approach
(cf. L.~E.~Dickson's proof that $SL_2(q)$ is two-generated for odd $q\ne 9$
\cite{Dickson}, or \cite{Doro}).


\section{Generators}

\noindent It is inconvenient to work with the quotient $\mathcal M/Z(\mathcal
M)$. We shall readily identify the elements $M$ and $-M$ in all of our
computations, and multiply by $-1$ freely. For $M\in\mathcal M$, a matrix, let
$M'$ denote the transpose of $M$.

In order to linearize our notation, we introduce two mappings $u$, $l:{\mathbb
F}^3\to\mathcal M$ defined by
\begin{displaymath}
    u(\alpha)=\uelm{\alpha},\myspace l(\alpha)=\lelm{\alpha}.
\end{displaymath}
Next consider $t:{\mathbb F}^3\setminus\{0\}\to{\mathbb F}^3$ given by
\begin{displaymath}
    t(\alpha_1,\,\alpha_2,\,\alpha_3)=\left\{\begin{array}{ll}
        (-\alpha_1^{-1},\, 0,\, 0)  & \textrm{if $\alpha_1\ne 0$,}\\
        (0,\, -\alpha_2^{-1},\, 0)  & \textrm{if $\alpha_1=0$, $\alpha_2\ne 0$,}\\
        (0,\, 0,\, -\alpha_3^{-1})  & \textrm{otherwise.}
    \end{array} \right.
\end{displaymath}
Note that $\alpha\cdot t(\alpha)=-1$. Finally, let $s(\alpha)$ stand for the
matrix
\begin{displaymath}
    \aelm{\alpha}{t(\alpha)}.
\end{displaymath}
Recall that, for matrices in $\mathcal M$,
\begin{displaymath}
    \melm{a}{\alpha}{\beta}{b}^{-1}=\melm{b}{-\alpha}{-\beta}{a},
\end{displaymath}
and also the two special cases of $(4.1)$, $(4.2)$ \cite{Paige}:
\begin{eqnarray}
    l(\alpha)&=&s(\alpha)'u(t(\alpha))(-s(\alpha)'),\label{Eq:l}\\
    u(\alpha)&=&s(\alpha)l(t(\alpha))(-s(\alpha)).\label{Eq:u}
\end{eqnarray}
Observe that $l(\alpha)^{-1}=l(-\alpha)$, and $u(\alpha)^{-1}=u(-\alpha)$.

We start our search for generators with the following result due to Paige:
\begin{proposition}\label{Pr:PaigeGen}
Every simple Moufang loop $\mathcal M/Z(\mathcal M)$ is generated by
\begin{displaymath}
    \{u(\alpha),\,l(\alpha);\; 0\ne\alpha\in{\mathbb F}^3\}.
\end{displaymath}
\end{proposition}
\begin{proof}
Combine Lemmas $4.2$ and $4.3$ \cite{Paige}.
\end{proof}

Let us identify all nonzero elements of $\mathbb{PF}^3={\mathbb F}^3/\mathbb F$
with those vectors from ${\mathbb F}^3$ whose first nonzero coordinate equals
1. Also, let $e_1=(1$, $0$, $0)$, $e_2=(0$, $1$, $0)$, and $e_3=(0$, $0$, $1)$,
as usual.
\begin{proposition}\label{Pr:PaigeProjGen}
Assume that $\mathbb F$ is a prime field. Then $\mathcal M/Z(\mathcal M)$ is
generated by
\begin{displaymath}
    \{u(\alpha),\, l(\alpha);\; \alpha\in\mathbb{PF}^3\}.
\end{displaymath}
\end{proposition}
\begin{proof}
First check that $u(a\alpha)u(b\alpha)=u((a+b)\alpha)$, and
$l(a\alpha)l(b\alpha)=l((a+b)\alpha)$ for all $\alpha\in {\mathbb F}^3$, $a$,
$b\in \mathbb F$.  Given $0\ne\beta\in{\mathbb F}^3$, there is $a\in\mathbb F$
and $\alpha\in \mathbb{PF}^3$ such that $\beta = a\alpha$. Since $\mathbb F$
is prime, we can use $a$ as an exponent, and write $u(\beta)=u(\alpha)^a$,
$l(\beta)=l(\alpha)^a$. We are finished, by Proposition \ref{Pr:PaigeGen}.
\end{proof}
\begin{proposition}\label{Pr:FirstReduction}
Assume that $\mathbb F$ is prime. Then $\mathcal M/Z(\mathcal M)$ is generated
by
\begin{displaymath}
    \mathcal S=\{u(e_1),\, u(e_2),\, u(e_3)\}\cup \{s(\alpha),\, s(\alpha)';\;
        \alpha\in\mathbb{PF}^3\}.
\end{displaymath}
\end{proposition}
\begin{proof}
Observe that, for $\alpha\in \mathbb{PF}^3$, we have $t(\alpha)\in\{-e_1$,
$-e_2$, $-e_3\}$. Given $l(\alpha)$ with $\alpha\in \mathbb{PF}^3$,
(\ref{Eq:l}) yields $l(\alpha)=s(\alpha)'u(-e_i)(-s(\alpha)')$ for some $e_i$,
$1\le i\le 3$. In particular, the elements $l(e_i)$, $1\le i\le 3$, are
generated by $\mathcal S$.

Symmetrically, given $u(\alpha)$ with $\alpha\in \mathbb{PF}^3$, equation
(\ref{Eq:u}) yields $u(\alpha)=s(\alpha)l(-e_i)(-s(\alpha))$ for some $e_i$,
$1\le i\le 3$; and we are done by Proposition \ref{Pr:PaigeProjGen}.
\end{proof}
\begin{lemma}\label{Lm:SameT}
Assume $\alpha$, $\beta\in{\mathbb F}^3\setminus\{0\}$ are such that
$t(\alpha)=t(\beta)$. Then $s(\beta)=s(\alpha)l(-\alpha\times\beta)$,
and $s(\beta)'=s(\alpha)'u(\alpha\times\beta)$.
\end{lemma}
\begin{proof}
By definition, we have
\begin{displaymath}
    s(\alpha)s(\beta)
    =\melm{\alpha\cdot t(\beta)}{-t(\alpha)\times t(\beta)}
        {\alpha\times\beta}{t(\alpha)\cdot\beta}
    =\melm{-1}{0}{\alpha\times\beta}{-1}.
\end{displaymath}
Thus $s(\beta)=-s(\alpha)^{-1}l(-\alpha\times\beta)
=s(\alpha)l(-\alpha\times\beta)$. As for the remaining equation, start with
$s(\alpha)'s(\beta)'=-u(\alpha\times\beta)$.
\end{proof}
\begin{proposition}\label{Pr:Big}
Assume $\mathbb F$ is prime. Then $\mathcal M/Z(\mathcal M)$ is generated by
\begin{displaymath}
    \{u(e_i),\, s(e_i);\; 1\le i\le 3\}.
\end{displaymath}
\end{proposition}
\begin{proof}
Thanks to Proposition \ref{Pr:FirstReduction}, we only need to generate
elements $s(\alpha)$ and $s(\alpha)'$, for $\alpha\in\mathbb{PF}^3$. First of
all, note that $s(e_i)'=-s(e_i)$, and, by $(\ref{Eq:l})$,
$l(e_i)=s(e_i)'u(e_i)^{-1}s(e_i)$, for $1\le i\le 3$. Also, $s(0,\,1,\,
a)=s(e_2)l(e_1)^{-a}$ and $s(0$, $1$, $a)'=s(e_2)'u(e_1)^a$ for all
$a\in\mathbb F$. Similarly, $s(1$, $a$, $0)=s(e_1)l(e_3)^{-a}$ and $s(1$, $a$,
$0)'=s(e_1)'u(e_3)^a$ for all $a\in\mathbb F$. Next, by $(\ref{Eq:l})$ and
$(\ref{Eq:u})$,
\begin{eqnarray*}
    l(1,a,0)&=&s(1,a,0)'u(-e_1)(-s(1,a,0)'),\\
    u(1,a,0)&=&s(1,a,0)l(-e_1)(-s(1,a,0)).
\end{eqnarray*}
For $0\ne a\in\mathbb F$, we have $l(a^{-1}$, $1$, $0)=l(1$, $a$,
$0)^{a^{-1}}$, $u(a^{-1}$, $1$, $0)=u(1,\,a,\,0)^{a^{-1}}$. Hence we obtained
$l(-a$, $1$, $0)$ and $u(-a$, $1$, $0)$ for all $a\in\mathbb F$. Finally,
Lemma \ref{Lm:SameT} yields
\begin{eqnarray*}
    s(1,a,b)&=&s(1,a,0)l(b(-a,1,0))=s(1,a,0)l(-a,1,0)^b,\\
    s(1,a,b)'&=&s(1,a,0)'u(b(a,-1,0))=s(1,a,0)'u(-a,1,0)^{-b},
\end{eqnarray*}
for all $b\in\mathbb F$.
\end{proof}

We may further reduce the number of generators one by one down to three. The
equations (\ref{Eq:Reductions}) below have been carefully chosen. They are
satisfied independently of $\char{\mathbb F}$, and are as simple as the author
was able to find. We leave the somewhat lengthy verification to the reader. Let
\begin{displaymath}
    x=\melm{0}{e_3}{-e_3}{1}.
\end{displaymath}
Then
\begin{equation}\label{Eq:Reductions}
\begin{array}{rcl}
    s(e_1)&=&s(e_3)s(e_2),\\
    s(e_2)&=&[u(e_1)u(e_3)\cdot u(e_2)u(e_1)][u(e_3)^{-1}\cdot u(e_2)s(e_3)],\\
    u(e_3)&=&x^{-1}s(e_3), \\
    s(e_3)&=&[u(e_2)u(e_1)\cdot x]
        [xu(e_1)]\cdot[x^2u(e_2)\cdot u(e_1)x^2u(e_1)].
\end{array}
\end{equation}
\begin{theorem}\label{Th:Main}
Assume that $\mathbb F$ is a prime field. Then the simple Moufang loop
$\mathcal M/Z(\mathcal M)$ is generated by the three elements
\begin{displaymath}
    \uelm{(1, 0, 0)},\myspace \uelm{(0, 1, 0)},
    \myspace \melm{0}{(0,0,1)}{(0,0,-1)}{1}.
\end{displaymath}
\end{theorem}
\section{Integral Cayley Numbers of Unit Norm}
\noindent The smallest nonassociative simple Moufang loop $\mathcal M_{120}$
of $120$ elements is constructed over the binary field. As Paige showed in
\cite{Paige}, it is isomorphic to the integral Cayley numbers of unit norm
modulo their center. The aim of this short section is to offer one explicit
isomorphism, and to conclude that the loop of integral Cayley numbers of unit
norm is generated multiplicatively by three elements. The isomorphism will
allow us to perform calculations inside integral Cayley numbers more
efficiently than by the conventional rules.

Unlike the case of complex numbers $\mathbb C$ and quaternions $\mathbb H$,
there are still many names for the eight-dimensional real algebra $\mathbb
O=\mathbb H\times\mathbb H$ obtained from $\mathbb H$ by the Cayley-Dickson
process: Cayley numbers, algebra of octaves, octonions. We prefer to use the
name \textit{octonions}. Recall that for $(q$, $Q)$, $(r,\,R)\in \mathbb
H\times\mathbb H$ the multiplication in $\mathbb O$ is defined by
\begin{displaymath}
    (q,\, Q)\cdot(r,\, R)=(qr-\overline{R}Q,\, Rq+Q\overline{r}),
\end{displaymath}
where the bar indicates conjugation in $\mathbb H$.
An alternative way to describe the multiplication is to introduce a new unit
$e$, regard $\mathbb O$ as $\mathbb H+\mathbb He$, and write
$(q+Qe)\cdot(r+Re)=(qr-\overline{R}Q)+(Rq+Q\overline{r})e$.
(For an excellent discussion concerning the notation in $\mathbb O$, see
\cite{Coxeter}.)

When $\mathbb O$ is viewed as an eight-dimensional real vector space, we will
use Dickson's notation for its basis, namely $1$, $i$, $j$, $k$, $e$, $ie$,
$je$, and $ke$. This is not the best choice when one wishes to describe the
multiplicative relations between basis elements in a compact way (see
\cite{Atlas}, \cite{Coxeter}), but it seems to be the best choice for what
follows.

The conjugate of the octonion $a=q+Qe$ is defined as $\overline{a}
=\overline{q}-Qe$. Its norm $N(a)$ is then the non-negative real number
$a\overline{a}$. For any $\mathbb F= \mathbb R$, $\mathbb C$, $\mathbb H$,
$\mathbb O$, a \textit{set of integral elements} of $\mathbb F$ is defined as
a maximal subset of $\mathbb F$ containing $1$, closed under multiplication
and subtraction, and such that both $N(a)$ and $a+\overline{a}$ are integers
for each element $a$ of the set. Such a set is unique for $\mathbb F=\mathbb
R$, $\mathbb C$, and $\mathbb H$. In the case of the octonions there are seven
such sets, all isomorphic. For the rest of this paper, we select the one which
Coxeter denotes by $J$, and calls \textit{integral Cayley numbers}.

By $J'$ we mean the $240$ elements of $J$ with norm one. See \cite{Coxeter},
p.$29$, for the list of all elements of $J'$. No three of the basis elements
generate $J'$. Following common practice, let $h=1/2\cdot(i+j+k+e)\in J'$.
Coxeter knew that $i$, $j$, and $h$ generated $J$ by multiplication and
subtraction, but he did not notice that $J'$ is generated by the above elements
just under multiplication. From \ref{Th:Main}, we know that there must be three
elements generating $J'/\{1,\,-1\}$. Indeed, $i$, $j$, and $h$ do the job. One
of the possible isomorphisms $\varphi:J'/\{1,\,-1\}\to\mathcal M_{120}$ is
determined by
\begin{displaymath}
   i\mapsto\aelm{e_3}{e_3},\myspace j\mapsto\aelm{e_2}{e_2},\myspace
         h\mapsto\melm{1}{(0,1,0)}{(1,0,1)}{1},
\end{displaymath}
which, however, is rather tedious to check by hand. The author acknowledges
that he used his own set of GAP $4$ libraries \cite{GAP} to confirm the
computation. Since $i^2=-1$, we see that $i$, $j$, and $h$ generate $J'$.

In order to be able to carry out calculations in $J'$, we would like to know
$\varphi(e)$.
\begin{lemma}
$(i)$ $e=-(jh\cdot hi)\cdot kh$ in $\mathbb O$.\newline
$(ii)$ The element corresponding to $e$ under $\varphi$ is
\begin{displaymath}
    \melm{0}{(1,1,1)}{(1,1,1)}{0}.
\end{displaymath}
\end{lemma}
\begin{proof}
Consider the multiplicative relations in \cite{Coxeter}, p.$567$ (or p.$28$).
In particular, we have $hi=-1-ih$. Hence
\begin{eqnarray*}
    -(jh\cdot hi)\cdot kh   &=& (jh+jh\cdot ih)\cdot kh = (jh+k-h-ih)\cdot kh \\
                    &=& jh\cdot kh+k\cdot kh-h\cdot kh-ih\cdot kh \\
                &=& (-i+h-kh)+(-h)-(k-h)-(j-h-kh)\\
                &=& -i-j-k+2h,
\end{eqnarray*}
which equals $e$, since $h=1/2\cdot(i+j+k+e)$.
\end{proof}

\end{article}

\begin{references}

\bibitem{Atlas} J.~H.~Conway and al., \textit{$\mathbb{ATLAS}$ of Finite Groups},
Oxford University Press, New York, 1985.

\bibitem{Bruck} R.~H.~Bruck, \textit{A Survey of Binary Systems}, Ergebnisse
der Mathematik und ihrer Grenzgebiete, Neue Folge, Heft $20$, Springer-Verlag,
1958.

\bibitem{Coxeter} H.~S.~M.~Coxeter, \textit{Integral Cayley Numbers}, Duke
Mathematical Journal, Vol. 13, No. 4, December, 1946. Reprinted in
H.~S.~M.~Coxeter, \textit{Twelve Geometric Essays}, Southern Illinois
University Press, 1968.

\bibitem{Dickson} L.~E.~Dickson, \textit{Linear Groups with an Exposition of
the Galois Field Theory}, Teubner, 1901; reprinted by Dover, 1958.

\bibitem{Doro} S.~Doro, \textit{Simple Moufang loops}, Math. Proc. Camb. Phil.
Soc. (1978), \textbf{83}, pp.377--392

\bibitem{GAP} The GAP Group, GAP --- Groups, Algorithms, and Programming,
Version 4.1; Aachen, St Andrews, 1999.
(http://www-gap.dcs.st-and.ac.uk/\~{}gap). The specific libraries related to
this paper are available at http://www.public.iastate.edu/\~{}petr.

\bibitem{Liebeck} M.~W.~Liebeck, \textit{The classification of finite simple
Moufang loops}, Math. Proc. Camb. Phil. Soc. (1987), {\bf 102}, pp.33--47.

\bibitem{Moufang} R.~Moufang, \textit{Zur Struktur von Alternativk\"orpern},
Math. Ann. \textbf{110} (1935), pp. 416--430.

\bibitem{Paige} L.~Paige, \textit{A Class of Simple Moufang Loops}, Proceedings
of the American Mathematical Society, Vol. $7$, Issue $3$, pp. $471$--$482$,
June, $1956$.

\bibitem{Pflugfelder} H.~O.~Pflugfelder, \textit{Quasigroups and Loops:
Introduction}, (Sigma series in pure mathematics; 7), Heldermann Verlag
Berlin, 1990.

\bibitem{Wilson} J.~S.~Wilson, \textit{Economical Generating Sets for Finite
Simple Groups}, London Mathematical Society Lecture Note Series \textbf{207},
Groups of Lie Type and Their Geometries, Cambridge University Press, 1995.

\end{references}
\end{document}